\documentclass[10pt]{article}

\usepackage{amsmath}
\usepackage{amsfonts}
\usepackage{amssymb}

\usepackage{xcolor}
\usepackage{graphicx}
\definecolor{LINK}{rgb}{0.35, 0.15, 0.1}
\definecolor{URL}{rgb}{0.6, 0.2, 0.7}

\usepackage{colortbl}
\definecolor{GR}{rgb}{0.0, 0.5, 0.0}
\definecolor{BR}{rgb}{0.7, 0.4, 0.0}
\definecolor{BL}{rgb}{0.0, 0.3, 1.0}
\definecolor{LI}{rgb}{0.6, 0.0, 0.8}
\definecolor{OR}{rgb}{0.7, 0.4, 0.0}
\definecolor{RO}{rgb}{0.8, 0.1, 0.3}

\newcommand{\TC}{\textcolor}

\usepackage{listings}
\usepackage[linkcolor=LINK, urlcolor=URL, citecolor=URL]{hyperref}
\hypersetup{linktocpage,colorlinks=true,pdfborder={0 0 0}}
\newcommand{\OEIS}[1]{\href{https://oeis.org/#1}{#1}}

\usepackage{geometry}
\title{New Results on the Stopping Time Behavior of the Collatz 3x + 1 Function}

\author{Mike Winkler}

\date{\small{\small Fakult\"at f\"ur Mathematik,\\Ruhr-Universit\"at Bochum, Germany,\\ mike.winkler@ruhr-uni-bochum.de\\www.mikematics.de\\[5mm]}}

\begin{document}
  
  \maketitle
  
  \begin{abstract}
	\noindent For the Collatz 3x + 1 function exists with $\sigma_n=\lfloor1+n\cdot\log_23\rfloor$ for each $n\in\mathbb{N}$ a set of different residue classes $(\text{mod}\ 2^{\sigma_n})$ of starting numbers $s$ with finite stopping time $\sigma(s)=\sigma_n$.
	Let $z_n$ be the number of these residue classes for each $n\geq0$ as listed in the OEIS as A100982.
	It is conjectured that for each $n\geq4$ the value of $z_n$ is given by
	\begin{align*}
	  z_n=\binom{m+n-2}{m}-\sum_{i=2}^{n-1}\binom{\big\lfloor\frac{3(n-i)+\delta}{2}\big\rfloor}{n-i}\cdot z_i,
	\end{align*}
	where $m=\big\lfloor(n-1)\cdot\log_23\big\rfloor-(n-1)$ and $\delta\in\mathbb{Z}$ assumes different values within the sum at intervals of 5 or 6 terms. This allows us to create an iterative algorithm which generates $z_n$ for each $n>6$.
	\\[70mm]
  \end{abstract}
  
  \newpage
  \tableofcontents
  
  \newpage
  
  \section{The Collatz $3x + 1$ conjecture}
  
  The Collatz $3x + 1$ function is defined as a function $T:\mathbb{N}\rightarrow\mathbb{N}$ on
  the set of positive integers by
  \begin{align*}
    T(x):=\left\{\begin{array}{lcr}T_0:=\displaystyle{\frac{x}{2}} & \mbox{if $x$ is even},\\ \\T_1:=\displaystyle{\frac{3x+1}{2}} & \mbox{if $x$ is odd}.\end{array}\right.
  \end{align*}
  \\ \\
  Let $T^0(s)=s$ and $T^k(s)=T\left(T^{k-1}(s)\right)$ for $k\in\mathbb{N}$. Then the Collatz sequence for $s\in\mathbb{N}$ is $C(s)=\left(T^k(s)\mid k=0,1,2,3,\dotsc\right)$.
  \\ \\
  For example, the starting number $s=11$ generates the Collatz sequence
    \[C(11)=(11,17,26,13,20,10,5,8,4,2,1,2,1,2,1,\dotsc).\]
  A Collatz sequence can only assume two possible forms. Either it falls into a cycle or it grows to infinity. The unproved conjecture to this problem is that each Collatz sequence enters the trivial cycle $(1,2,1,2,\dotsc)$.
  \\ \\
  \textit{Note}: Proofs for the Theorems shown in this article can be found in \textsc{Winkler}\cite{Winkler2017}.
  \\
  
  \section{Stopping time}
  
  \subsection{The stopping time $\sigma(s)$}
  
  Collatz's conjecture is equivalent to the conjecture that for each $s\in\mathbb{N},s>1$, there exists $k\in\mathbb{N}$ such that $T^k(s)<s$. The least $k\in\mathbb{N}$ such that $T^k(s)<s$ is called the stopping time of $s$, which we will denote by $\sigma(s)$. It is not hard to verify that
  \\ \\
  \hspace*{10mm} $\sigma(s)=1$ \quad if \, $s\equiv0\ (\text{mod}\ 2)$,\\
  \hspace*{10mm} $\sigma(s)=2$ \quad if \, $s\equiv1\ (\text{mod}\ 4)$,\\
  \hspace*{10mm} $\sigma(s)=4$ \quad if \, $s\equiv3\ (\text{mod}\ 16)$,\\
  \hspace*{10mm} $\sigma(s)=5$ \quad if \, $s\equiv11,23\ (\text{mod}\ 32)$,\\
  \hspace*{10mm} $\sigma(s)=7$ \quad if \, $s\equiv7,15,59\ (\text{mod}\ 128)$,\\
  \hspace*{10mm} $\sigma(s)=8$ \quad if \, $s\equiv39,79,95,123,175,199,219\ (\text{mod}\ 256)$,\\
  \hspace*{10mm} and so forth.
  \\ \\
  Let $\sigma_n=\lfloor1+n\cdot\log_23\rfloor$ then generally applies for all $n\in\mathbb{N},\,n\geq0$, that
  \begin{align*}
    \sigma(s)=\sigma_n \quad \text{if} \quad s\equiv x_1,x_2,x_3,\dotsc,x_z\ (\text{mod}\ 2^{\sigma_n}).
  \end{align*}
  Let $z_n$ be the number of residue classes $(\text{mod}\ 2^{\sigma_n})$ for each $n\geq0$. Then we have
  \begin{align*}
    z_0=1,\ z_1=1,\ z_2=1,\ z_3=2,\ z_4=3,\ z_5=7,\ z_6=12,\ z_7=30,\ \dotsc
  \end{align*}
  \\
  \noindent\textbf{Theorem 1}: Let $\sigma_n=\lfloor1+n\cdot\log_23\rfloor$ then for each $n\in\mathbb{N},\,n\geq0$, there exists a set of $z_n$ residue classes $(\text{mod}\ 2^{\sigma_n})$ of starting numbers $s$ with finite stopping time $\sigma(s)=\sigma_n$. \textit{(End of theorem)}
  \\
  The possible stopping times $\sigma(s)$ are listed in \OEIS{A020914}. The associated residue classes $(\text{mod}\ 2^{\sigma_n})$ are listed in \OEIS{A177789}. The number of residue classes $z_n$ for $n\geq1$ are listed in \OEIS{A100982}.
  \\ \\
  Appendix 5.4 shows a list of the first residue classes $(\text{mod}\ 2^{\sigma_n})$ up to $\sigma(s)=15$.
  \\
  
  \subsection{A stopping time term formula for odd $s$}
  
  \textbf{Theorem 2}: Let $C^a(s)=\left(T^k(s)\mid k=0,\dotsc,a\right)$ with $a\geq1$ be a finite subsequence of $C(s)$, and let $\sigma_n=\lfloor1+n\cdot\log_23\rfloor$. Then for each $n\in\mathbb{N}$ an odd starting number $s$ has the stopping time $\sigma(s)=\sigma_n$, if the appropriate subsequence $C^{\sigma_n-1}(s)$ consists of $n$ odd terms, and $\alpha_i=k$, if and only if $T^k(s$) in $C^{\sigma_n-1}(s)$ is odd. Then it is
  \begin{align}
    T^{\sigma_n}(s)=\frac{3^n}{2^{\sigma_n}}\cdot s+\sum_{i=1}^n \frac{3^{n-i}2^{\alpha_i}}{2^{\sigma_n}}<T^0(s).
  \end{align}
  \textit{(End of theorem)}
  \\ \\
  \textit{Example}: For $n=4$ there is $\sigma_4=\lfloor1+4\cdot\log_23\rfloor=7$. For $s=59$ we get by equation (1)
  \begin{align*}
    T^{7}(59)=\frac{3^4}{2^7}\cdot59+ \frac{3^3 2^0+3^2 2^1+3^1 2^3+3^0 2^4}{2^7}=38<59.
  \end{align*}
  \\
  \textit{Explanation}: The subsequence $C^6(59)=(59,89,134,67,101,152,76)$ consists of  four $(n=4)$ odd terms $59,89,67,101$. The powers of two $\alpha_i$ yield as follows: $T^0=59$ is odd, so $\alpha_1=0$. $T^1=89$ is odd, so $\alpha_2=1$. $T^2=134$ is even. $T^3=67$ is odd, so $\alpha_3=3$. $T^4=101$ is odd, so $\alpha_4=4$. $T^5=152$ is even. $T^6=76$ is even.
  \\ \\
  \textit{Note}: There is $\sigma(s)=7$ not only for $s=59$, but for every $s\equiv59\ (\text{mod}\ 2^7)$.
  \\
  
  \section{Diophantine equations}
  
  \subsection{Subsequences and their binary simplification}
  
  \textbf{Theorem 3}: Let $m=\big\lfloor(n-1)\cdot\log_23\big\rfloor-(n-1)$ and $\sigma_n=\lfloor1+n\cdot\log_23\rfloor$ for each $n\in\mathbb{N},\,n\geq4$. According to the conditions of theorem 2, if an odd starting number $s$ has the stopping time $\sigma(s)=\sigma_n$, then the first $m+n$ terms in $C(s)$ represents sufficiently the stopping time of $s$, because all further terms are even before the term $T^{\sigma_n}(s)<s$ is reached. \textit{(End of theorem)}
  \\ \\
  To simplify the distribution of the even and odd terms in $C(s)$ let "0" represents an even term and "1" represents an odd term.
  \\ \\
  \textit{Example}: For $n=4$ there is $m=\big\lfloor3\cdot\log_23\big\rfloor-3=1$ and $\sigma_4=\lfloor1+4\cdot\log_23\rfloor=7$. Then for the subsequences $C^{\sigma_n}(s)$ there is
  \begin{align*}
    C^7(7)&=(7,11,17,26,13,20,10,5) & \text{simplified by} \quad (1,1,1,0,1,0,0,1),\\
    C^7(15)&=(15,23,35,53,80,40,20,10) & \text{simplified by} \quad (1,1,1,1,0,0,0,0),\\
    C^7(59)&=(59,89,134,67,101,152,76,38) & \text{simplified by} \quad (1,1,0,1,1,0,0,0).
  \end{align*}
  And for the sufficiently subsequences $C^{m+n-1}(s)$ there is
  \begin{align*}
    C^4(7)&=(7,11,17,26,13) & \text{simplified by} \quad (1,1,1,0,1),\\
    C^4(15)&=(15,23,35,53,80) & \text{simplified by} \quad (1,1,1,1,0),\\
    C^4(59)&=(59,89,134,67,101) & \text{simplified by} \quad (1,1,0,1,1).
  \end{align*}
  
  \subsection{Binary tuples}
  
  Let $m=\big\lfloor(n-1)\cdot\log_23\big\rfloor-(n-1)$ for each $n\in\mathbb{N},\,n\geq4$.
  \\ \\
  Let $\mathcal{A}(n)$ be a binary $(m+n-2)$-tuple for each $n\in\mathbb{N},\,n\geq4$, defined by
    \[\mathcal{A}(n)=(a_1, \cdots, a_m,a_{m+1}, \cdots, a_{m+n-2}),\]
    \[with \quad a_1, \cdots, a_m:=0 \quad and \quad a_{m+1}, \cdots, a_{m+n-2}:=1.\]
  \\
  Let $\mathbb{A}_j(n)$ be the set of $\mathcal{A}(n)$ and all permutations in lexicographic ordering of $\mathcal{A}(n)$ for each $n\in\mathbb{N},\,n\geq4$, whereby $j$ is the number of all these tuples for each $n$, calculated by 
    \[j=\binom{m+n-2}{m}=\frac{(m+n-2)!}{m!\cdot(n-2)!},\]\\
  which generates the sequence listed in \OEIS{A293308}.
  \\ \\
  Let $\mathcal{B}(n)$ be a binary $(m+n)$-tuple for each $n\in\mathbb{N},\,n\geq4$, defined by
    \[\mathcal{B}(n)=(b_1, \cdots, b_{m+n}) \quad with \quad b_1:=1, \ b_2:=1.\]
  Let $\mathbb{B}_j(n)$ be the set of $j$ tuples $\mathcal{B}(n)$ for each $n\in\mathbb{N},\,n\geq4$, where $b_3, \cdots, b_{m+n}$ is equal to exactly one tuple of $\mathbb{A}_j(n)$.
  \\ \\
  \textit{Example}: For $n=7$ there is $m=3$ and $j=56$. Therefore we get the 8-tuple $\mathcal{A}(7)=(0,0,0,1,1,1,1,1)$. There are further 55 permutations in lexicographic ordering of $\mathcal{A}(7)$. Because of this $\mathbb{A}_{56}(7)$ contains 56 different 8-tuples, and $\mathbb{B}_{56}(7)$ contains 56 different 10-tuples  $\mathcal{B}(7)$. Appendix 5.2 shows these example in detail.
  \\
  
  \subsection{The Diophantine equations and their solutions}
  
  According to theorem 2, the behavior of a Collatz sequence is clearly related to the way in which the powers of 2 are distributed among the powers of 3 in the term $3^{n-i}2^{\alpha_i}$ of equation (1).
  \\ \\
  According to theorem 3, for each $n\in\mathbb{N},n\geq4$, only for the $j$ binary tuples of $\mathbb{B}_j(n)$ the conditions of theorem 2 and equation (1) are complied.
  \\ \\
  \textbf{Theorem 4}:  Let $\sigma_n=\lfloor1+n\cdot\log_23\rfloor$. By interpreting the binary tuples of $\mathbb{B}_j(n)$ as such a binary simplification for the even and odd terms in $C^{m+n-1}(s)$, there exists for each $n\in\mathbb{N},\,n\geq4$, according to theorem 2, for each binary tuple of $\mathbb{B}_j(n)$ a Diophantine equation 
  \begin{align}
    y=\frac{3^n}{2^{\sigma_n}}\cdot x+\sum_{i=1}^n \frac{3^{n-i}2^{\alpha_i}}{2^{\sigma_n}},
  \end{align}
  which has exactly one integer solution $(x,y)$ for $0<x<2^{\sigma_n}$, so that $y=T^{\sigma_n}(x)<x$ in $C(x)$. \textit{(End of theorem)}
  \\ \\
  \textbf{Theorem 5}: For each $n\geq4$ the $j$ solutions $(x,y)$ represent $j$ different residue classes of starting numbers $s$ with stopping time $4\leq\sigma(s)\leq\sigma_n$. More precise: For each $n\geq4$ the $j$ solutions represent $z_n$ residue classes of starting numbers $s$ with the stopping time $\sigma(s)=\sigma_n$, and $j-z_n$ residue classes of starting numbers $s$ with stopping time $4\leq\sigma(s)<\sigma_n$. For the $z_n$ solutions $(x,y)$ with $\sigma(x)=\sigma_n$, there is $x$ the smallest number of the residue class $[x]_{2^{\sigma_n}}$. \textit{(End of theorem)}
  \\ \\
  \textit{Example}: For $n=5$ there are $j=10$ integer solutions $(x,y)$. The next table shows the 10 binary tuples of $\mathbb{B}_{10}(5)$ and the appropriate integer solution $(x,y)$.
  \\
  \small
  \lstset{numbers=none}
  \begin{lstlisting}
             (1, 1, 0, 0, 1, 1, 1)   (211, 202)
             (1, 1, 0, 1, 0, 1, 1)   (107, 103)
             (1, 1, 0, 1, 1, 0, 1)   (123, 118)
             (1, 1, 0, 1, 1, 1, 0)   (219, 209)
             (1, 1, 1, 0, 0, 1, 1)   (183, 175)
             (1, 1, 1, 0, 1, 0, 1)   (199, 190)
             (1, 1, 1, 0, 1, 1, 0)   ( 39,  38)
             (1, 1, 1, 1, 0, 0, 1)   ( 79,  76)
             (1, 1, 1, 1, 0, 1, 0)   (175, 167)
             (1, 1, 1, 1, 1, 0, 0)   ( 95,  91)
  \end{lstlisting}
  \normalsize
  \quad\\
  Of these 10 solutions, there are $z_5=7$ solutions for which $x$ has the stopping time $\sigma(x)=\sigma_5=8$, and $10-7=3$ solutions for which $x$ has stopping time $4\leq\sigma(x)<8$.
  \\ \\
  It is $\ 211\in[3]_{16}$, $\ 107\in[11]_{32}$ and $\ 183\in[23]_{32}$. So these 10 solutions represent 10 different residue classes of starting numbers $s$ with stopping time $4\leq\sigma(s)\leq8$. According to Chapter 2, there is
  \\ \\
  \hspace*{10mm} $\sigma(s)=4$ \quad if \, $s\equiv3\ (\text{mod}\ 16)$,\\
  \hspace*{10mm} $\sigma(s)=5$ \quad if \, $s\equiv11,23\ (\text{mod}\ 32)$,\\
  \hspace*{10mm} $\sigma(s)=8$ \quad if \, $s\equiv39,79,95,123,175,199,219\ (\text{mod}\ 256)$.
  \\ \\
  Note that for $\sigma(s)=8$ the 7 residue classes are equal to the values of $x$.
  \\ \\
  \\
  The next table shows the distribution of the solutions $(x,y)$ to the possible stopping times $\sigma(x)$ for $n=4,\dots,12$.
  \small
  \begin{center}
    \begin{tabular}{|c|c|r|r|r|r|r|r|r|r|r|}
      \hline & $n$&  4 & 5 &  6 &  7 &  8 &   9&  10 &   11 &   12\\
      \hline $\sigma(x)$&&&&    &    &    &    &     &      &     \\
      \hline  4 & &  0 & 1 &  1 &  6 & 28 &  36& 165 &  220 & 1001\\
      \hline  5 & &  0 & 2 &  2 & 10 & 42 &  56& 240 &  330 & 1430\\
      \hline  7 & &  3 & 0 &  0 &  3 & 15 &  18&  84 &  108 &  495\\
      \hline  8 & &    & 7 &  0 &  7 & 28 &  35& 147 &  196 &  840\\
      \hline 10 & &    &   & 12 &  0 & 12 &  12&  60 &   72 &  336\\
      \hline 12 & &    &   &    & 30 &  0 &   0&  30 &   30 &  180\\
      \hline 13 & &    &   &    &    & 85 &   0&  85 &   85 &  425\\
      \hline 15 & &    &   &    &    &    & 173&   0 &    0 &  173\\
      \hline 16 & &    &   &    &    &    &    & 476 &    0 &  476\\
      \hline 18 & &    &   &    &    &    &    &     &  961 &    0\\
      \hline 20 & &    &   &    &    &    &    &     &      & 2652\\
      \hline\hline Sum&& 3& 10& 15& 56& 210& 330& 1287& 2002& 8008\\
      \hline
    \end{tabular}
  \end{center}
  \normalsize
  \textit{How to read}: "Sum" is equal to the value of $j$ and means the sum of the values of $\sigma(x)$ for each $n$. For $n=5$ there are 1+2+7=10 integer solutions $(x,y)$. From these solutions, there are 7 with $\sigma(x)=8$, 2 with $\sigma(x)=5$ and 1 with $\sigma(x)=4$. No entry is equal to "0".
  \\ \\
  The PARI/GP programs 1 - 3 in Appendix 5.1 make the algorithm of this chapter more clearly. Program 2 generates the table of page 6, and program 3 computes the values of the upper table. Appendix 5.3 shows $\mathbb{B}_j(n)$ with integer solutions $(x,y)$ for $n=4,\dotsc,7$.
  \\
  
  \subsection{Analysis of the solutions}
  
  An analysis of the table entries shows that for each $\sigma(x)$ the number of solutions for all $n\geq4$ are integer multiples of
  \TC{RO}{$z_n$}.
  \small
  \begin{center}
    \begin{tabular}{|c|c|r|r|r|r|r|r|r|r|r|}
      \hline&$n$&4&5&6&7&8&9&10&11&12\\
      \hline$\sigma(x)$&&&&&&&&&&\\
      \hline4&&$0\cdot1$&$1\cdot1$&$1\cdot1$&$\TC{LI}6\cdot1$&$\TC{BR}{28}\cdot1$&$\TC{LI}{36}\cdot1$&$\TC{BR}{165}\cdot1$&$\TC{LI}{220}\cdot1$&$\TC{BR}{1001}\cdot1$\\
      \hline5&&$0\cdot2$&$1\cdot2$&$1\cdot2$&$\TC{LI}5\cdot2$&$\TC{BR}{21}\cdot2$&$\TC{LI}{28}\cdot2$&$\TC{BR}{120}\cdot2$&$\TC{LI}{165}\cdot2$&$\TC{BR}{715}\cdot2$\\
      \hline7&&$1\cdot\TC{RO}3$&$0\cdot3$&$0\cdot3$&$\TC{LI}1\cdot3$&$\TC{BR}5\cdot3$&$\TC{LI}6\cdot3$&$\TC{BR}{28}\cdot3$&$\TC{LI}{36}\cdot3$&$\TC{BR}{165}\cdot3$\\
      \hline8&&&$1\cdot\TC{RO}7$&$0\cdot7$&$\TC{LI}1\cdot7$&$\TC{BR}4\cdot7$&$\TC{LI}5\cdot7$&$\TC{BR}{21}\cdot7$&$\TC{LI}{28}\cdot7$&$\TC{BR}{120}\cdot7$\\
      \hline10&&&&$1\cdot\TC{RO}{12}$&$\TC{LI}0\cdot12$&$\TC{BR}1\cdot12$&$\TC{LI}1\cdot12$&$\TC{BR}5\cdot12$&$\TC{LI}{6}\cdot12$&$\TC{BR}{28}\cdot12$\\
      \hline1&&&&&$1\cdot\TC{RO}{30}$&$\TC{LI}0\cdot30$&$0\cdot30$&$\TC{LI}1\cdot30$&$\TC{BR}1\cdot30$&$\TC{LI}6\cdot30$\\
      \hline13&&&&&&$1\cdot\TC{RO}{85}$&$0\cdot85$&$\TC{LI}1\cdot85$&$\TC{BR}1\cdot85$&$\TC{LI}5\cdot85$\\
      \hline15&&&&&&&$1\cdot\TC{RO}{173}$&$\TC{LI}0\cdot173$&$0\cdot173$&$\TC{LI}1\cdot173$\\
      \hline1&&&&&&&&$1\cdot\TC{RO}{476}$& $0\cdot476$&$\TC{LI}1\cdot476$\\
      \hline18&&&&&&&&&$1\cdot\TC{RO}{961}$&$\TC{LI}0\cdot961$\\
      \hline20&&&&&&&&&&$1\cdot\TC{RO}{2652}$\\
      \hline\hline Sum&&3&10&15&56&210&330&1287&2002&8008\\
      \hline
    \end{tabular}
  \end{center}
  \normalsize
  
  \noindent\\ \\It is not hard to recognize that the colored multiplying factors of $z_n$ are numbers of Pascal's triangle. It seems that these numbers are given for each $n$ by the following binomial coefficients for $k\in\mathbb{N}$.
  
  \begin{align*}
    \binom{\big\lfloor\frac{3k-2}{2}\big\rfloor}{k}=& \ \TC{LI}{0,1,1,5,6,28,36,165,220,1001,1365,6188,8568,38760,\dotsc}
    \\ \\
    \binom{\big\lfloor\frac{3k-1}{2}\big\rfloor}{k}=& \ \TC{BR}{1,1,4,5,21,28,120,165,715,1001,4368,6188,27132,\dotsc}
    \\ \\
    \binom{\big\lfloor\frac{3k}{2}\big\rfloor}{k}=& \ \TC{BL}{1,3,4,15,21,84,120,495,715,3003,4368,18564,27132,\dotsc}
  \end{align*}
  \\ \\
  \noindent Now we conjecture that for each $n\geq4$ the value of $z_n$ is equal to the difference of $j$ and the sum of $n-2$ terms, where each term is given by the product of a special binomial coefficient and a number $z_n$ for $n=2,\dotsc,n-1$. Let
  \begin{align}
    z_n=\binom{m+n-2}{m}-\sum_{i=2}^{n-1}\binom{\big\lfloor\frac{3(n-i)+\delta}{2}\big\rfloor}{n-i}\cdot z_i,
  \end{align}
  where $m=\big\lfloor(n-1)\cdot\log_23\big\rfloor-(n-1)$ and $\delta\in\mathbb{Z}$ assumes different values within the sum at intervals of 5 or 6 terms. The example for $n=13$ on page 9 and the algorithm on page 10 will make this more clear.
  \\ \\
  The next table shows the conjectured values for $n=13,\dotsc,18$. These values are not computed with program 3, but generated with program 4.
  \small
  \begin{center}
    \begin{tabular}{|c|c|r|r|r|r|r|r|}
      \hline&$n$&13&14&15&16&17&18\\
      \hline$\sigma(x)$&&&&&&&\\
      \hline4&&$\TC{BL}{4368}\cdot1$&$\TC{BR}{6188}\cdot1$&$\TC{BL}{27132}\cdot1$&$\TC{BR}{38760}\cdot1$&$\TC{BL}{170544}\cdot1$&$\TC{BR}{245157}\cdot1$\\
      \hline5&&$\TC{BL}{3003}\cdot2$&$\TC{BR}{4368}\cdot2$&$\TC{BL}{18564}\cdot2$&$\TC{BR}{27132}\cdot2$&$\TC{BL}{116280}\cdot2$&$\TC{BR}{170544}\cdot2$\\
      \hline7&&$\TC{BL}{715}\cdot3$&$\TC{BR}{1001}\cdot3$&$\TC{BL}{4368}\cdot3$&$\TC{BR}{6188}\cdot3$&$\TC{BL}{27132}\cdot3$&$\TC{BR}{38760}\cdot3$\\
      \hline8&&$\TC{BL}{495}\cdot7$&$\TC{BR}{715}\cdot7$&$\TC{BL}{3003}\cdot7$&$\TC{BR}{4368}\cdot7$ &$\TC{BL}{18564}\cdot7$&$\TC{BR}{27132}\cdot7$\\
      \hline10&&$\TC{BL}{120}\cdot12$&$\TC{BR}{165}\cdot12$&$\TC{BL}{715}\cdot12$&$\TC{BR}{1001}\cdot12$&$\TC{BL}{4368}\cdot12$&$\TC{BR}{6188}\cdot12$\\
      \hline12&&$\TC{BR}{28}\cdot30$&$\TC{LI}{36}\cdot30$&$\TC{BR}{165}\cdot30$&$\TC{LI}{220}\cdot30$&$\TC{BR}{1001}\cdot30$&$\TC{LI}{1365}\cdot30$\\
      \hline13&&$\TC{BR}{21}\cdot85$&$\TC{LI}{28}\cdot85$&$\TC{BR}{120}\cdot85$&$\TC{LI}{165}\cdot85$&$\TC{BR}{715}\cdot85$&$\TC{LI}{1001}\cdot85$\\
      \hline15&&$\TC{BR}{5}\cdot173$&$\TC{LI}{6}\cdot173$&$\TC{BR}{28}\cdot173$&$\TC{LI}{36}\cdot173$&$\TC{BR}{165}\cdot173$&$\TC{LI}{220}\cdot173$\\
      \hline16&&$\TC{BR}{4}\cdot476$&$\TC{LI}{5}\cdot476$&$\TC{BR}{21}\cdot476$&$\TC{LI}{28}\cdot476$&$\TC{BR}{120}\cdot476$&$\TC{LI}{165}\cdot476$\\
      \hline18&&$\TC{BR}{1}\cdot961$&$\TC{LI}{1}\cdot961$&$\TC{BR}{5}\cdot961$&$\TC{LI}{6}\cdot961$&$\TC{BR}{28}\cdot961$&$\TC{LI}{36}\cdot961$\\
      \hline20&& $0\cdot2652$&$0\cdot2652$&$\TC{LI}{1}\cdot2652$&${1}\cdot2652$&$\TC{LI}{6}\cdot2652$&${7}\cdot2652$\\
      \hline21&&$1\cdot\TC{RO}{8045}$&${0}\cdot8045$&$\TC{LI}{1}\cdot8045$&${1}\cdot8045$&$\TC{LI}{5}\cdot8045$&${6}\cdot8045$\\
      \hline23&&&$1\cdot\TC{RO}{17637}$&$0\cdot17637$&$0\cdot17637$&$\TC{LI}{1}\cdot17637$&${1}\cdot17637$\\
      \hline24&&&&$1\cdot\TC{RO}{51033}$&$0\cdot51033$&$\TC{LI}{1}\cdot51033$&$1\cdot51033$\\
      \hline26&&&&&$1\cdot\TC{RO}{108950}$&$0\cdot108950$&$0\cdot108950$\\
      \hline27&&&&&&$1\cdot\TC{RO}{312455}$&$0\cdot312455$\\
      \hline29&&&&&&&$1\cdot\TC{RO}{663535}$\\
      \hline\hline Sum&&31824&50388&203490&319770&1307504&2042975\\
      \hline
    \end{tabular}
  \end{center}
  \normalsize\quad
  \\ \\
  \noindent\textit{Example}: This example shows the working of equation (3) for $n=13$.
  \\
  \begin{align*}
    z_n&=\binom{m+n-2}{m}-\sum_{i=2}^{n-1}\binom{\big\lfloor\frac{3(n-i)+\delta}{2}\big\rfloor}{n-i}\cdot z_i \quad\text{with}\quad m=\big\lfloor(n-1)\cdot\log_23\big\rfloor-(n-1)
    \\ \\
    \TC{RO}{z_{n}}&=\frac{(m+n-2)!}{m!\cdot(n-2)!}-\TC{BL}{\sum_{i=2}^{6}\binom{\big\lfloor\frac{3(n-i)+0}{2}\big\rfloor}{n-i}}\cdot z_i-\TC{BR}{\sum_{i=7}^{11}\binom{\big\lfloor\frac{3(n-i)-1}{2}\big\rfloor}{n-i}}\cdot z_i-\sum_{i=12}^{12}\binom{\big\lfloor\frac{3(n-i)-2}{2}\big\rfloor}{n-i}\cdot z_i
    \\ \\
    \TC{RO}{z_{13}}&=\frac{(7+13-2)!}{7!\cdot(13-2)!}-\TC{BL}{\sum_{i=2}^{6}\binom{\big\lfloor\frac{3(13-i)}{2}\big\rfloor}{13-i}}\cdot z_i-\TC{BR}{\sum_{i=7}^{11}\binom{\big\lfloor\frac{3(13-i)-1}{2}\big\rfloor}{13-i}}\cdot z_i-\binom{\big\lfloor\frac{3(13-12)-2}{2}\big\rfloor}{13-12}\cdot z_{12}
    \\ \\
    \TC{RO}{z_{13}}&=\frac{18!}{7!\cdot 11!}-\TC{BL}{\sum_{i=2}^{6}\binom{\big\lfloor\frac{39-3i}{2}\big\rfloor}{13-i}}\cdot z_i-\TC{BR}{\sum_{i=7}^{11}\binom{\big\lfloor\frac{38-3i}{2}\big\rfloor}{13-i}}\cdot z_i-\binom{\big\lfloor\frac{1}{2}\big\rfloor}{1}\cdot z_{12}
    \\ \\
    \TC{RO}{z_{13}}&=31824-\TC{BL}{4368}\cdot 1-\TC{BL}{3003}\cdot 2-\TC{BL}{715}\cdot 3-\TC{BL}{495}\cdot 7-\TC{BL}{120}\cdot 12-\TC{BR}{28}\cdot 30\\&\quad\quad\quad\ \ \, -\TC{BR}{21}\cdot 85-\TC{BR}{5}\cdot 173-\TC{BR}{4}\cdot 476-\TC{BR}{1}\cdot961-0\cdot2652
    \\ \\
    \TC{RO}{z_{13}}&=31824-17424-6355=\TC{RO}{8045}.
  \end{align*}
  
  \subsection{An iterative stopping time algorithm}
  
  The results of the analysis of the solutions and equation (3) enables us to devise an iterative algorithm which generates with only five initial numbers $z_2,\dotsc,z_6$ each further number $z_n$ and also the number of all integer solutions $(x,y)$ for $0<x<2^{\sigma_n}$ from the Diophantine equation (2) with same stopping time $\sigma(x)$ as seen in the tables on page 7 and 8.
  \\ \\
  The next PARI/GP program 4 shows this algorithm, which outputs a list of the numbers $z_n$ for $6<n\leq limit$, as listed in the OEIS as a \href{https://oeis.org/A100982/list}{simple list} of A100982.
  \\ \\
  The correctness of this algorithm for the values of $z_n$ has been proved with another algorithm\footnote{see Theorem 2 in \OEIS{A100982}} for each $n<54$.
  \\
  \lstset{numbers=left}
  \begin{lstlisting}
  {/* Program 4 - stopping time algorithm */
   limit=54;
   zn=vector(limit);
   /* 5 initial numbers */
   zn[2]=1; zn[3]=2; zn[4]=3; zn[5]=7; zn[6]=12;
   /* main algorithm */
   f=1; e1=-1; e2=-2; 
   for(n=7, limit,
     m=floor((n-1)*log(3)/log(2))-(n-1);
     j=(m+n-2)!/(m!*(n-2)!);
     if(n>6*f, if(frac(n/2)==0, e=e1, e=e2));
     if(frac((n-6 )/12)==0, f++; e1=e1+2);
     if(frac((n-12)/12)==0, f++; e2=e2+2);
     Sum=a=b=0; c=1; d=5;
     until(c>=n-1,
       for(i=2+a*5+b, 1+d+a*5,
         if(i>11 && frac((i+2)/6)==0, b++);
         delta=e-a;
         Sum=Sum+binomial(floor((3*(n-i)+delta)/2),n-i)*zn[i];
         c++;
       );
       a++;
       for(k=3, 50, if(n>=k*6 && a==k-1, d=k+3));
     );
   zn[n]=j-Sum;
   print(n" "zn[n]);
   );
  }
  \end{lstlisting}
  
  \noindent\\ The lines 10 and 19 of program 4 show the main calculations of equation (3). The only complicated thing is the computation of the pattern of the lower and upper bounds of the index of summation $i$ and the appropriate value of $\delta$.
  \\ \\
  The upper limit of the \textit{for}-loop in line 23 depends on the value of $n$ or \textit{limit} and must be large enough. The bigger the value of these upper limit, the longer the runtime of the algorithm.
  \\ \\
  The PARI/GP program 5 in Appendix 5.1 shows this algorithm with a different output as seen in the tables on page 7 and 8. For other algorithms which generates the values of $z_n$ by a different way see \OEIS{A100982} and \textsc{Winkler}\cite{Winkler2017}.
  \\
  
  \section{Conclusion and some answers}
  
  \subsection{A general stopping time theorem}
  
  We see that the fact, if an odd starting number $s$ has finite stopping time is only dependent on how the first $m$ even and $n$ odd terms are distributed in $C^{m+n-1}(s)$ and $C^{\sigma_n-1}(s)$. Now we are able to formulate a general stopping time theorem without the use of equation (1).
  \\ \\
  \textbf{Theorem 6}: Let $m=\big\lfloor(n-1)\cdot\log_23\big\rfloor-(n-1)$ and $\sigma_n=\lfloor1+n\cdot\log_23\rfloor$. For all $n\in\mathbb{N},\,n\geq4$, an odd starting number $s$ has stopping time $4\leq\sigma(s)\leq\sigma_n$, if the binary simplification of the even and odd terms in the subsequence $C^{m+n-1}(s)$ is equal to a binary tuple of $\mathbb{B}_j(n)$ \textit{and} the subsequence $C^{\sigma_n-1}(s)$ consists of $n$ odd terms. \textit{(End of theorem)}
  \\
  
  \subsection{\small Why the Collatz conjecture is a Diophantine equation problem?}
  
  According to theorem 1, the Collatz conjecture is true, if the set of the residue classes $(\text{mod}\ 2^{\sigma_n})$ of starting numbers for all $n\geq0$ is equal to $\mathbb{N}$.  According to theorem 4 and 5, all possible residue classes for the stopping times are given by the Diophantine equations as their integer solutions. Therefore the Collatz conjecture is true, if there exists for each necessary residue class an appropriate Diophantine equation.
  \\
  
  \subsection{\small What are the consequences, if there exists another cycle than (1,2) or a sequence with infinite growth?}
  
  Then, according to theorem 6, there exist an odd starting number $s$ without stopping time. This starting number $s$ has the property that for each $n\in\mathbb{N},\,n\geq4$, there is no accordance of the binary simplification of the even and odd terms in the subsequence $C^{m+n-1}(s)$ and a binary tuple of $\mathbb{B}_j(n)$, or the subsequence $C^{\sigma_n-1}(s)$ consists not of $n$ odd terms.
  \\
  
  \subsection{\small Why has a small starting number like $s=27$ a so comparatively large stopping time $\sigma(27)=59$?}
  
  Because out of the finite sets of Diophantine equations for each $n\geq4$ not until $n=37$ a Diophantine equation has the solution $x=27$ or the residue class $27\ (\text{mod}\ 2^{59})$. Another answer: Not until $n=37$ the binary simplification of the even and odd terms in the subsequence $C^{m+n-1}(27)$ is equal to a binary tuple of $\mathbb{B}_j(n)$ while the subsequence $C^{\sigma_n-1}(27)$ consists of $n$ odd terms.
  \\
  
  \subsection{Future work on the Collatz graph}
  
  My future work is about the self-reference of the Collatz graph. I have devised two different directed graphs, which illustrate that there exists a self-reference of odd subsequences $C^a(s)$, clarified by the same background colour between both graphs and arrows in the first graph.
  \\ \\
  This means on the one hand that odd subsequences $C^a(s)$ with a small number of steps to reach 1 are connected with subsequences $C^a(2s\pm1)$ or $C^a(4s\pm3)$ with a big number of steps to reach 1. For example $C^3(15)=(15,23,35,53)$ and $C^3(31)=(31,47,71,107)$.
  \\ \\
  And on the other hand shows the first graph why and where to find numbers in the second graph. For example shows the green arrow in the first graph why $C^1(27)=(27,41)$ must hit the subsequence $C^4(31)=(31,47,71,107,161)$. In the second graph the subsequence $C^6(27)=(27,41,31,47,71,107,161)$ is located in column 9.
  \\ \\
  This is only a very brief and insufficient exposition of the theory of self-reference. For details of the construction of the second graph and its correctness see \textsc{Winkler}\cite{Winkler2010}. See this \href{http://mikematics.de/collatz-tree1.pdf}{link}\footnote{http://mikematics.de/collatz-tree1.pdf} for a picture of the second graph which contains all odd numbers up to 341, and this \href{http://mikematics.de/collatz-tree2.pdf}{link}\footnote{http://mikematics.de/collatz-tree2.pdf} for the complete structure up to column 6.
  \\ \\
  \begin{figure}[ht]
  	\centering
    \includegraphics[scale=0.5]{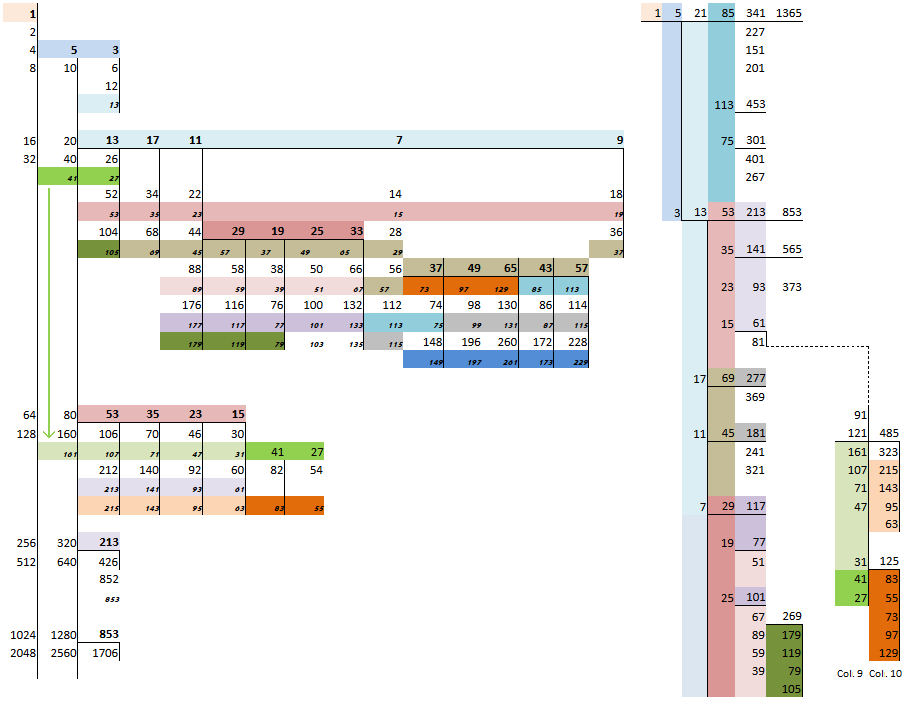}
    \caption{Two directed Collatz graphs}
  \end{figure}
  
  \newpage
  
  \section{Appendix}
  
  \subsection{Algorithms in PARI/GP}
  
  The programs 1 - 3 use the function "NextPermutation(a)", which generates all permutations in lexicographic ordering of a binary tuple $\mathcal{A}(n)$.
  
  \begin{lstlisting}
  NextPermutation(a)=
  {
   i=#a-1;
   while(!(i<1 || a[i]<a[i+1]), i--);
   if(i<1, return(0));
   k=#a;
   while(!(a[k]>a[i]), k--);
   t=a[k]; a[k]=a[i]; a[i]=t;
   for(k=i+1, (#a+i)/2,
     t=a[k]; a[k]=a[#a+1+i-k]; a[#a+1+i-k]=t;
   );
   return(a); 
  }
  \end{lstlisting}
  
  \noindent\\\\Program 1 computes all integer solutions $(x,y)$ for $0<x<2^{\sigma_n}$ from the Diophantine equation (2) for $n=7$.
  \begin{lstlisting}
  {/* Program 1 */
   n=7;
   m=floor((n-1)*log(3)/log(2))-(n-1);
   Sigma=floor(1+n*log(3)/log(2));
   j=(m+(n-2))!/(m!*(n-2)!);
   /* vectors for binary tuple */
   Alpha=[0,0,0,0,0,0,0];
       A=[0,0,0,1,1,1,1,1];
       B=[1,1,0,0,0,0,0,0,0,0];
   z=0; q=0;
   until(z==j, z++;
     /* generate B(n) from A(n) */
     for(i=1, m+n-2, B[2+i]=A[i]);
     /* determine the n values for Alpha[i] */
     i=1; for(k=1, m+n, if(B[k]==1, Alpha[i]=k-1; i++));
     /* calculate integer solutions from Diophantine equation */
     for(x=0, 2^Sigma,
       Sum=0; for(i=1, n, Sum=Sum+3^(n-i)*2^Alpha[i]);
       y=(3^n*x+Sum)/2^Sigma;
       if(frac(y)==0, print(B"   ("x", "y")"); q++);
     );
     A=NextPermutation(A);
   );
   print("There are "q" integer solutions (x,y) for n=7.");
  }
  \end{lstlisting}
  
  \newpage
  
  \noindent Program 2 computes all integer solutions $(x,y)$ for $0<x<2^{\sigma_n}$ from the Diophantine equation (2) for $2\leq n\leq limit$.
  \begin{lstlisting}
  {/* Program 2 */
   limit=100;
   for(n=2, limit,
     m=floor((n-1)*log(3)/log(2))-(n-1);
     Sigma=floor(1+n*log(3)/log(2));
     j=(m+(n-2))!/(m!*(n-2)!);
     /* generate vectors for binary tuple */
     Alpha=[]; for(i=1, n, Alpha=concat(Alpha,i));
               for(i=1, n, Alpha[i]=0);
     A=[]; for(i=1, m+n-2, A=concat(A,i));
           for(i=1, m+n-2, if(i<=m, A[i]=0, A[i]=1));
     B=[]; for(i=1, m+n, B=concat(B,i));
           for(i=1, m+n, if(i<=2, B[i]=1, B[i]=0));
     z=0; q=0;
     until(z==j, z++;
       /* generate B(n) from A(n) */
       for(i=1, m+n-2, B[2+i]=A[i]);
       /* determine the n values for Alpha[i] */
       i=1; for(k=1, m+n, if(B[k]==1, Alpha[i]=k-1; i++));
       /* calculate integer solutions from Diophantine equation */
       for(x=0, 2^Sigma,
         Sum=0; for(i=1, n, Sum=Sum+3^(n-i)*2^Alpha[i]);
         y=(3^n*x+Sum)/2^Sigma;
         if(frac(y)==0, print(B"   ("x", "y")"); q++);
       );
       A=NextPermutation(A);
     );
     print("There are "q" integer solutions (x,y) for n="n".");
     print;
   );
  }
  \end{lstlisting}
  
  \noindent\\Program 3 computes all integer solutions $(x,y)$ for $0<x<2^{\sigma_n}$ from the Diophantine equation (2) for $2\leq n\leq10$ and counts the solutions with same stopping time $\sigma(x)$.
  \begin{lstlisting}
  {/* Program 3 */
   for(n=2, 10,
     m=floor((n-1)*log(3)/log(2))-(n-1);
     Sigma=floor(1+n*log(3)/log(2));
     j=(m+(n-2))!/(m!*(n-2)!);
     /* initialize stopping time counters */
     w1=w2=w3=w4=w5=w6=w7=w8=w9=w10=0;
     /* generate vectors for binary tuple */
     Alpha=[]; for(i=1, n, Alpha=concat(Alpha,i));
               for(i=1, n, Alpha[i]=0);
     A=[]; for(i=1, m+n-2, A=concat(A,i));
           for(i=1, m+n-2, if(i<=m, A[i]=0, A[i]=1));
     B=[]; for(i=1, m+n, B=concat(B,i));
           for(i=1, m+n, if(i<=2, B[i]=1, B[i]=0));
     z=0; q=0;
     until(z==j, z++;
       /* generate B(n) from A(n) */
       for(i=1, m+n-2, B[2+i]=A[i]);
       /* determine the n values for Alpha[i] */
       i=1; for(k=1, m+n, if(B[k]==1, Alpha[i]=k-1; i++));
       /* calculate integer solutions from Diophantine equation */
       for(x=0, 2^Sigma,
         Sum=0; for(i=1, n, Sum=Sum+3^(n-i)*2^Alpha[i]);
         y=(3^n*x+Sum)/2^Sigma;
         if(frac(y)==0, q++;
           /* calculate stopping time ST of x */
           ST=0; c=x;
           until(c<x, if(frac(c/2)<>0, ST++; c=(3*c+1)/2);
             if(frac(c/2)==0, ST++; c=c/2;
               if((frac(c/2)<>0) && (c<x), break);
             );
           );   
           if(ST==2, w1++); if(ST==4, w2++); if(ST==5, w3++);
           if(ST==7, w4++); if(ST==8, w5++); if(ST==10, w6++);
           if(ST==12, w7++); if(ST==13, w8++); if(ST==15, w9++);
           if(ST==16, w10++);
         );
       );
       A=NextPermutation(A);
     );
     print;
     print("There are "q" integer solutions (x,y) for n="n".");
     print("number of x with ST= 2: "w1);
     print("number of x with ST= 4: "w2);
     print("number of x with ST= 5: "w3);
     print("number of x with ST= 7: "w4);
     print("number of x with ST= 8: "w5);
     print("number of x with ST=10: "w6);
     print("number of x with ST=12: "w7);
     print("number of x with ST=13: "w8);
     print("number of x with ST=15: "w9);
     print("number of x with ST=16: "w10);
   );
  }
  \end{lstlisting}
  \normalsize
  Appendix 5.3 shows the output of program 2 for $n=4,\dotsc,7$.
  \\ \\
  \noindent Program 5 is equal to program 4 in chapter 3.5 but outputs the number of all integer solutions in a form as seen in the tables on page 7 and 8.
  \begin{lstlisting}
  {/* Program 5 - stopping time algorithm */
   limit=54;
   zn=vector(limit);
   /* 5 initial numbers */
   zn[2]=1; zn[3]=2; zn[4]=3; zn[5]=7; zn[6]=12;
   /* main algorithm */
   f=1; e1=-1; e2=-2; 
   for(n=7, limit,
     m=floor((n-1)*log(3)/log(2))-(n-1);
     j=(m+n-2)!/(m!*(n-2)!);
     if(n>6*f, if(frac(n/2)==0, e=e1, e=e2));
     if(frac((n-6 )/12)==0, f++; e1=e1+2);
     if(frac((n-12)/12)==0, f++; e2=e2+2);
     Sum=a=b=0; c=1; d=5;
     until(c>=n-1,
       for(i=2+a*5+b, 1+d+a*5,
         if(i>11 && frac((i+2)/6)==0, b++);
         delta=e-a;
         s=binomial(floor((3*(n-i)+delta)/2),n-i);
         Sum=Sum+s*zn[i];
         c++;
         if(zn[i]!=0, print(i"   "s"  "zn[i]"   "s*zn[i]));
         );
       a++;
       for(k=3, 50, if(n>=k*6 && a==k-1, d=k+3));
     );
     zn[n]=j-Sum;
     print(n" "zn[n]);
   );
  }
  \end{lstlisting}
  
  \newpage
  
  \subsection{Permutations in lexicographic ordering for $n=7$}
  
  \normalsize For $n=7$ there is $\mathcal{A}(7)=(0,0,0,1,1,1,1,1)$,
  \\ \\
  $\mathbb{A}_{56}(7)=$
  \lstset{basicstyle=\scriptsize\ttfamily,numbers=none}
  \begin{lstlisting}
			      { (0,0,0,1,1,1,1,1) ,     (1,0,1,0,1,1,0,1) ,
				      (0,0,1,0,1,1,1,1) ,     (1,0,1,0,1,1,1,0) ,
				      (0,0,1,1,0,1,1,1) ,     (1,0,1,1,0,0,1,1) ,
				      (0,0,1,1,1,0,1,1) ,     (1,0,1,1,0,1,0,1) ,
				      (0,0,1,1,1,1,0,1) ,     (1,0,1,1,0,1,1,0) ,
				      (0,0,1,1,1,1,1,0) ,     (1,0,1,1,1,0,0,1) ,
				      (0,1,0,0,1,1,1,1) ,     (1,0,1,1,1,0,1,0) ,
				      (0,1,0,1,0,1,1,1) ,     (1,0,1,1,1,1,0,0) ,
				      (0,1,0,1,1,0,1,1) ,     (1,1,0,0,0,1,1,1) ,
				      (0,1,0,1,1,1,0,1) ,     (1,1,0,0,1,0,1,1) ,
				      (0,1,0,1,1,1,1,0) ,     (1,1,0,0,1,1,0,1) ,
				      (0,1,1,0,0,1,1,1) ,     (1,1,0,0,1,1,1,0) ,
				      (0,1,1,0,1,0,1,1) ,     (1,1,0,1,0,0,1,1) ,
				      (0,1,1,0,1,1,0,1) ,     (1,1,0,1,0,1,0,1) ,
				      (0,1,1,0,1,1,1,0) ,     (1,1,0,1,0,1,1,0) ,
				      (0,1,1,1,0,0,1,1) ,     (1,1,0,1,1,0,0,1) ,
				      (0,1,1,1,0,1,0,1) ,     (1,1,0,1,1,0,1,0) ,
				      (0,1,1,1,0,1,1,0) ,     (1,1,0,1,1,1,0,0) ,
				      (0,1,1,1,1,0,0,1) ,     (1,1,1,0,0,0,1,1) ,
				      (0,1,1,1,1,0,1,0) ,     (1,1,1,0,0,1,0,1) ,
				      (0,1,1,1,1,1,0,0) ,     (1,1,1,0,0,1,1,0) ,
				      (1,0,0,0,1,1,1,1) ,     (1,1,1,0,1,0,0,1) ,
				      (1,0,0,1,0,1,1,1) ,     (1,1,1,0,1,0,1,0) ,
				      (1,0,0,1,1,0,1,1) ,     (1,1,1,0,1,1,0,0) ,
				      (1,0,0,1,1,1,0,1) ,     (1,1,1,1,0,0,0,1) ,
				      (1,0,0,1,1,1,1,0) ,     (1,1,1,1,0,0,1,0) ,
				      (1,0,1,0,0,1,1,1) ,     (1,1,1,1,0,1,0,0) ,
				      (1,0,1,0,1,0,1,1) ,     (1,1,1,1,1,0,0,0)  }
  \end{lstlisting}
  
  \noindent$\mathbb{B}_{56}(7)=$
  \begin{lstlisting}
			  { (1,1,0,0,0,1,1,1,1,1) , (1,1,1,0,1,0,1,1,0,1) ,
			    (1,1,0,0,1,0,1,1,1,1) , (1,1,1,0,1,0,1,1,1,0) ,
			    (1,1,0,0,1,1,0,1,1,1) , (1,1,1,0,1,1,0,0,1,1) ,
			    (1,1,0,0,1,1,1,0,1,1) , (1,1,1,0,1,1,0,1,0,1) ,
			    (1,1,0,0,1,1,1,1,0,1) , (1,1,1,0,1,1,0,1,1,0) ,
			    (1,1,0,0,1,1,1,1,1,0) , (1,1,1,0,1,1,1,0,0,1) ,
			    (1,1,0,1,0,0,1,1,1,1) , (1,1,1,0,1,1,1,0,1,0) ,
			    (1,1,0,1,0,1,0,1,1,1) , (1,1,1,0,1,1,1,1,0,0) ,
			    (1,1,0,1,0,1,1,0,1,1) , (1,1,1,1,0,0,0,1,1,1) ,
			    (1,1,0,1,0,1,1,1,0,1) , (1,1,1,1,0,0,1,0,1,1) ,
			    (1,1,0,1,0,1,1,1,1,0) , (1,1,1,1,0,0,1,1,0,1) ,
			    (1,1,0,1,1,0,0,1,1,1) , (1,1,1,1,0,0,1,1,1,0) ,
			    (1,1,0,1,1,0,1,0,1,1) , (1,1,1,1,0,1,0,0,1,1) ,
			    (1,1,0,1,1,0,1,1,0,1) , (1,1,1,1,0,1,0,1,0,1) ,
			    (1,1,0,1,1,0,1,1,1,0) , (1,1,1,1,0,1,0,1,1,0) ,
			    (1,1,0,1,1,1,0,0,1,1) , (1,1,1,1,0,1,1,0,0,1) ,
			    (1,1,0,1,1,1,0,1,0,1) , (1,1,1,1,0,1,1,0,1,0) ,
			    (1,1,0,1,1,1,0,1,1,0) , (1,1,1,1,0,1,1,1,0,0) ,
			    (1,1,0,1,1,1,1,0,0,1) , (1,1,1,1,1,0,0,0,1,1) ,
			    (1,1,0,1,1,1,1,0,1,0) , (1,1,1,1,1,0,0,1,0,1) ,
			    (1,1,0,1,1,1,1,1,0,0) , (1,1,1,1,1,0,0,1,1,0) ,
			    (1,1,1,0,0,0,1,1,1,1) , (1,1,1,1,1,0,1,0,0,1) ,
			    (1,1,1,0,0,1,0,1,1,1) , (1,1,1,1,1,0,1,0,1,0) ,
			    (1,1,1,0,0,1,1,0,1,1) , (1,1,1,1,1,0,1,1,0,0) ,
			    (1,1,1,0,0,1,1,1,0,1) , (1,1,1,1,1,1,0,0,0,1) ,
			    (1,1,1,0,0,1,1,1,1,0) , (1,1,1,1,1,1,0,0,1,0) ,
			    (1,1,1,0,1,0,0,1,1,1) , (1,1,1,1,1,1,0,1,0,0) ,
			    (1,1,1,0,1,0,1,0,1,1) , (1,1,1,1,1,1,1,0,0,0)  }
  \end{lstlisting}
  
  \newpage
  
  \subsection{$\mathbb{B}_j(n)$ with integer solutions $(x,y)$ for $n=4,\dotsc,7$}
  This is also the output of program 2 for $n=4,\dotsc,7$.
  \\
  \lstset{basicstyle=\scriptsize\ttfamily,
          keywordstyle=\bfseries\ttfamily}
  \begin{lstlisting}
    (1, 1, 0, 1, 1)   (59, 38)
    (1, 1, 1, 0, 1)   (7, 5)
    (1, 1, 1, 1, 0)   (15, 10)
    There are 3 integer solutions (x,y) for n=4.
  \end{lstlisting}
  \begin{lstlisting}
    (1, 1, 0, 0, 1, 1, 1)   (211, 202)
    (1, 1, 0, 1, 0, 1, 1)   (107, 103)
    (1, 1, 0, 1, 1, 0, 1)   (123, 118)
    (1, 1, 0, 1, 1, 1, 0)   (219, 209)
    (1, 1, 1, 0, 0, 1, 1)   (183, 175)
    (1, 1, 1, 0, 1, 0, 1)   (199, 190)
    (1, 1, 1, 0, 1, 1, 0)   (39, 38)
    (1, 1, 1, 1, 0, 0, 1)   (79, 76)
    (1, 1, 1, 1, 0, 1, 0)   (175, 167)
    (1, 1, 1, 1, 1, 0, 0)   (95, 91)
    There are 10 integer solutions (x,y) for n=5.
  \end{lstlisting}
  \begin{lstlisting}
    (1, 1, 0, 0, 1, 1, 1, 1)   (595, 425)
    (1, 1, 0, 1, 0, 1, 1, 1)   (747, 533)
    (1, 1, 0, 1, 1, 0, 1, 1)   (507, 362)
    (1, 1, 0, 1, 1, 1, 0, 1)   (347, 248)
    (1, 1, 0, 1, 1, 1, 1, 0)   (923, 658)
    (1, 1, 1, 0, 0, 1, 1, 1)   (823, 587)
    (1, 1, 1, 0, 1, 0, 1, 1)   (583, 416)
    (1, 1, 1, 0, 1, 1, 0, 1)   (423, 302)
    (1, 1, 1, 0, 1, 1, 1, 0)   (999, 712)
    (1, 1, 1, 1, 0, 0, 1, 1)   (975, 695)
    (1, 1, 1, 1, 0, 1, 0, 1)   (815, 581)
    (1, 1, 1, 1, 0, 1, 1, 0)   (367, 262)
    (1, 1, 1, 1, 1, 0, 0, 1)   (735, 524)
    (1, 1, 1, 1, 1, 0, 1, 0)   (287, 205)
    (1, 1, 1, 1, 1, 1, 0, 0)   (575, 410)
    There are 15 integer solutions (x,y) for n=6.
  \end{lstlisting}
  \begin{lstlisting}    
    (1, 1, 0, 0, 0, 1, 1, 1, 1, 1)   (3523, 1883)
    (1, 1, 0, 0, 1, 0, 1, 1, 1, 1)   (3827, 2045)
    (1, 1, 0, 0, 1, 1, 0, 1, 1, 1)   (1299, 695)
    (1, 1, 0, 0, 1, 1, 1, 0, 1, 1)   (979, 524)
    (1, 1, 0, 0, 1, 1, 1, 1, 0, 1)   (2131, 1139)
    (1, 1, 0, 0, 1, 1, 1, 1, 1, 0)   (2899, 1549)
    (1, 1, 0, 1, 0, 0, 1, 1, 1, 1)   (3979, 2126)
    (1, 1, 0, 1, 0, 1, 0, 1, 1, 1)   (1451, 776)
    (1, 1, 0, 1, 0, 1, 1, 0, 1, 1)   (1131, 605)
    (1, 1, 0, 1, 0, 1, 1, 1, 0, 1)   (2283, 1220)
    (1, 1, 0, 1, 0, 1, 1, 1, 1, 0)   (3051, 1630)
    (1, 1, 0, 1, 1, 0, 0, 1, 1, 1)   (187, 101)
    (1, 1, 0, 1, 1, 0, 1, 0, 1, 1)   (3963, 2117)
    (1, 1, 0, 1, 1, 0, 1, 1, 0, 1)   (1019, 545)
    (1, 1, 0, 1, 1, 0, 1, 1, 1, 0)   (1787, 955)
    (1, 1, 0, 1, 1, 1, 0, 0, 1, 1)   (1755, 938)
    (1, 1, 0, 1, 1, 1, 0, 1, 0, 1)   (2907, 1553)
    (1, 1, 0, 1, 1, 1, 0, 1, 1, 0)   (3675, 1963)
    (1, 1, 0, 1, 1, 1, 1, 0, 0, 1)   (1435, 767)
    (1, 1, 0, 1, 1, 1, 1, 0, 1, 0)   (2203, 1177)
    (1, 1, 0, 1, 1, 1, 1, 1, 0, 0)   (2587, 1382)
    (1, 1, 1, 0, 0, 0, 1, 1, 1, 1)   (2007, 1073)
    (1, 1, 1, 0, 0, 1, 0, 1, 1, 1)   (3575, 1910)
    (1, 1, 1, 0, 0, 1, 1, 0, 1, 1)   (3255, 1739)
    (1, 1, 1, 0, 0, 1, 1, 1, 0, 1)   (311, 167)
    (1, 1, 1, 0, 0, 1, 1, 1, 1, 0)   (1079, 577)
    (1, 1, 1, 0, 1, 0, 0, 1, 1, 1)   (2311, 1235)
    (1, 1, 1, 0, 1, 0, 1, 0, 1, 1)   (1991, 1064)
    (1, 1, 1, 0, 1, 0, 1, 1, 0, 1)   (3143, 1679)
    (1, 1, 1, 0, 1, 0, 1, 1, 1, 0)   (3911, 2089)
    (1, 1, 1, 0, 1, 1, 0, 0, 1, 1)   (3879, 2072)
    (1, 1, 1, 0, 1, 1, 0, 1, 0, 1)   (935, 500)
    (1, 1, 1, 0, 1, 1, 0, 1, 1, 0)   (1703, 910)
    (1, 1, 1, 0, 1, 1, 1, 0, 0, 1)   (3559, 1901)
    (1, 1, 1, 0, 1, 1, 1, 0, 1, 0)   (231, 124)
    (1, 1, 1, 0, 1, 1, 1, 1, 0, 0)   (615, 329)
    (1, 1, 1, 1, 0, 0, 0, 1, 1, 1)   (3727, 1991)
    (1, 1, 1, 1, 0, 0, 1, 0, 1, 1)   (3407, 1820)
    (1, 1, 1, 1, 0, 0, 1, 1, 0, 1)   (463, 248)
    (1, 1, 1, 1, 0, 0, 1, 1, 1, 0)   (1231, 658)
    (1, 1, 1, 1, 0, 1, 0, 0, 1, 1)   (1199, 641)
    (1, 1, 1, 1, 0, 1, 0, 1, 0, 1)   (2351, 1256)
    (1, 1, 1, 1, 0, 1, 0, 1, 1, 0)   (3119, 1666)
    (1, 1, 1, 1, 0, 1, 1, 0, 0, 1)   (879, 470)
    (1, 1, 1, 1, 0, 1, 1, 0, 1, 0)   (1647, 880)
    (1, 1, 1, 1, 0, 1, 1, 1, 0, 0)   (2031, 1085)
    (1, 1, 1, 1, 1, 0, 0, 0, 1, 1)   (2143, 1145)
    (1, 1, 1, 1, 1, 0, 0, 1, 0, 1)   (3295, 1760)
    (1, 1, 1, 1, 1, 0, 0, 1, 1, 0)   (4063, 2170)
    (1, 1, 1, 1, 1, 0, 1, 0, 0, 1)   (1823, 974)
    (1, 1, 1, 1, 1, 0, 1, 0, 1, 0)   (2591, 1384)
    (1, 1, 1, 1, 1, 0, 1, 1, 0, 0)   (2975, 1589)
    (1, 1, 1, 1, 1, 1, 0, 0, 0, 1)   (1087, 581)
    (1, 1, 1, 1, 1, 1, 0, 0, 1, 0)   (1855, 991)
    (1, 1, 1, 1, 1, 1, 0, 1, 0, 0)   (2239, 1196)
    (1, 1, 1, 1, 1, 1, 1, 0, 0, 0)   (383, 205)
    There are 56 integer solutions (x,y) for n=7.
  \end{lstlisting}
  
  \subsection{Stopping time residue classes up to $\sigma(s)=15$}
  
  \tiny\noindent
  $\sigma(s)=1$\\
  if \ $s\equiv$ \ 0 \ $(\text{mod}\ 2)$\\
  \\
  $\sigma(s)=2$\\
  if \ $s\equiv$ \ 1 \ $(\text{mod}\ 4)$\\
  \\
  $\sigma(s)=4$\\
  if \ $s\equiv$ \ 3 \ $(\text{mod}\ 16)$\\
  \\
  $\sigma(s)=5$\\
  if \ $s\equiv$ \ 11, 23 \ $(\text{mod}\ 32)$\\
  \\
  $\sigma(s)=7$\\
  if \ $s\equiv$ \ 7, 15, 59 \ $(\text{mod}\ 128)$\\
  \\
  $\sigma(s)=8$\\
  if \ $s\equiv$ \ 39, 79, 95, 123, 175, 199, 219 \ $(\text{mod}\ 256)$\\
  \\
  $\sigma(s)=10$\\
  if \ $s\equiv$ \ 287, 347, 367, 423, 507, 575, 583, 735, 815, 923, 975, 999 \ $(\text{mod}\ 1024)$\\
  \\
  $\sigma(s)=12$\\
  if \ $s\equiv$ \ 231, 383, 463, 615, 879, 935, 1019, 1087, 1231, 1435, 1647, 1703, 1787, 1823, 1855, 2031, 2203, 2239, 2351, 2587, 2591, 2907, 2975, 3119, 3143, 3295, 3559, 3675, 3911, 4063 \ $(\text{mod}\ 4096)$\\
  \\
  $\sigma(s)=13$\\
  if \ $s\equiv$ \ 191, 207, 255, 303, 539, 543, 623, 679, 719, 799, 1071, 1135, 1191, 1215, 1247, 1327, 1563, 1567, 1727, 1983, 2015, 2075, 2079, 2095, 2271, 2331, 2431, 2607, 2663, 3039, 3067, 3135, 3455, 3483, 3551, 3687, 3835, 3903, 3967, 4079, 4091, 4159, 4199, 4223, 4251, 4455, 4507, 4859, 4927, 4955, 5023, 5103, 5191, 5275, 5371, 5439, 5607, 5615, 5723, 5787, 5871, 5959, 5979, 6047, 6215, 6375, 6559, 6607, 6631, 6747, 6815, 6983, 7023, 7079, 7259, 7375, 7399, 7495, 7631, 7791, 7847, 7911, 7967, 8047, 8103 \ $(\text{mod}\ 8192)$\\
  \\
  $\sigma(s)=15$\\
  if \ $s\equiv$ \ 127, 411, 415, 831, 839, 1095, 1151, 1275, 1775, 1903, 2119, 2279, 2299, 2303, 2719, 2727, 2767, 2799, 2847, 2983, 3163, 3303, 3611, 3743, 4007, 4031, 4187, 4287, 4655, 5231, 5311, 5599, 5631, 6175, 6255, 6503, 6759, 6783, 6907, 7163, 7199, 7487, 7783, 8063, 8187, 8347, 8431, 8795, 9051, 9087, 9371, 9375, 9679, 9711, 9959, 10055, 10075, 10655, 10735, 10863, 11079, 11119, 11567, 11679, 11807, 11943, 11967, 12063, 12143, 12511, 12543, 12571, 12827, 12967, 13007, 13087, 13567, 13695, 13851, 14031, 14271, 14399, 14439, 14895, 15295, 15343, 15839, 15919, 16027, 16123, 16287, 16743, 16863, 16871, 17147, 17727, 17735, 17767, 18011, 18639, 18751, 18895, 19035, 19199, 19623, 19919, 20079, 20199, 20507, 20527, 20783, 20927, 21023, 21103, 21223, 21471, 21727, 21807, 22047, 22207, 22655, 22751, 22811, 22911, 22939, 23231, 23359, 23399, 23615, 23803, 23835, 23935, 24303, 24559, 24639, 24647, 24679, 25247, 25503, 25583, 25691, 25703, 25831, 26087, 26267, 26527, 26535, 27111, 27291, 27759, 27839, 27855, 27975, 28703, 28879, 28999, 29467, 29743, 29863, 30311, 30591, 30687, 30715, 30747, 30767, 30887, 31711, 31771, 31899, 32155, 32239, 32575, 32603 \ $(\text{mod}\ 32768)$\\
  \\
  and so forth.
  \\
  
  \newpage
  \normalsize
  
  \section{References}
  
  \begingroup
  \renewcommand{\section}[2]{}

\end{document}